# Fuzzy Rough Relations

T. K. Samanta and Biswajit Sarkar

Department of Mathematics, Uluberia College, India - 711315.

e-mail: mumpu_tapas5@yahoo.co.in

Department of Mathematics, Bantul Mahakali High School, India - 711312

e-mail: chotonsarkar@yahoo.co.in

**Abstract**

*In this paper, the definition of fuzzy rough relation on a set will be introduced and then it would be proved that the collection of such relations is closed under different binary compositions such as, algebraic sum, algebraic product etc. Also the definitions of reflexive, symmetric and transitive fuzzy rough relations on a set are given and a few properties of them will be investigated. Lastly, we define a operation $\circ$ , which is a composition of two fuzzy rough relations, with the help of maxmin relation and thereafter it is shown that the collection of such relations is closed under the operation $\circ$ .*

## 1 Introduction

Most of problems in engineering, medical science, economics, environments etc. have various uncertainties. To overcome these uncertainties, some kind of theories were given like theory of fuzzy set [4], intuitionistic fuzzy sets [6], rough sets [8], soft sets [1] etc., which we can use as a mathematical tools for dealing with uncertainties. In 1965, Zadeh [4] initiated the novel concept of fuzzy set theory, thereafter in 1982, the concept of rough set theory was first given by Pawlak [8] and then in 1999, Molodtsov [1] initiated the concept of soft theory, all these are used for modeling incomplete knowledge, vagueness and uncertainties. In fact, all these concepts having a good application in other disciplines and real life problems are now catching momentum. But, it is seen that all these theories have their own difficulties, that is why in this paper we are going to study fuzzy rough relation depending upon fuzzy rough set, which is another one new mathematical tool for dealing with uncertainties.

In fact, there have been also attempts to "fuzzify" various mathematical structures like topological spaces, groups, rings, etc. and also concepts like relations, measure, probability, and automata etc. In this paper we will develop the concepts of relations with the help of fuzzy rough set. The concept of fuzzy relation on a set was defined by Zadeh [4, 5] and several authors have considered it further, has generalized the concept by considering fuzzy relations on fuzzy sets and developed the study of fuzzy graphs, obtaining beautiful analogs of several graph theoretical properties. In 1971, fuzzy relations were applied to consider clustering analysis. In this paper, we will introduce the definition of fuzzy rough relation on a set and then it would be proved that the collection of such relations is closed under different binary compositions such as, algebraic sum, algebraic product etc. Also we define reflexive, symmetric and transitive fuzzy rough relations on a set and investigate a few properties of them. Lastly, we define a operation ∘ with the help of maxmin relation and thereafter it is shown that the collection of such relations is closed under the operation ∘.

## 2 Preliminaries

Let $A = (U, R)$ be approximation space. The product product space is also an approximation space $A^2 = (U^2, S)$, where the indiscernibility relation $S \subseteq U^2$ is defined by $((x_1, y_1), (x_2, y_2)) \in S$ if and only if $(x_1, x_2) \in R$ and $(y_1, y_2) \in R$, for $x_1, x_2, y_1, y_2 \in U$. It can be easily verified that $S$ is an equivalence relation on $U$. The elements $(x_1, y_1)$ and $(x_2, y_2)$ are indiscernible in $S$ if and only if the elements $x_1, x_2$ are indiscernible in $R$ and so are $y_1, y_2$ in $R$. This implies that the equivalence class containing the element $(x, y)$ with respect to $S$, denoted by $[x, y]_S$, should be equal to the Cartesian product of $[x]_R$ and $[y]_R$.

The concepts of rough set can be easily extended to a relation, mainly due to the fact that a relation is also a set, i.e. a subset of a Cartesian product. So, let $A = (U, R)$ be an approximation space. Let $X \subseteq U$. A relation $T$ on $X$ is said to be a rough relation on $X$ if $\underline{T} \neq \overline{T}$, where $\underline{T}$ and $\overline{T}$ are lower and upper approximation of $T$, respectively defined by

$$\underline{T} = \{ (x, y) \in U \times U : [x, y]_S \subseteq X \}$$

$$\overline{T} = \{ (x, y) \in U \times U : [x, y]_S \cap X \neq \varphi \}$$

**Definition 2.1** *A fuzzy subset $A$ of the universe $U$ is characterized by a membership function $\mu_A : U \to [0, 1]$. Then the product $A \times A$ is defined by the membership function $\mu_{A \times A}(x, y) = \min \{ \mu_A(x), \mu_A(y) \}$. A fuzzy relation $R$ on $A$ is a fuzzy subset of $A \times A$. i.e., $\mu_R(x, y) \leq \mu_{A \times A}(x, y) \quad \forall x, y \in U$.*

There are three kind of fuzzy relations mainly, which are as follows :

**1. Reflexive Fuzzy Relation** : A fuzzy relation $R$ on $A$ is said to be *reflexive* iff $\mu_R(x, x) = 1 \ \forall x \in U$ such that $\mu_A(x) > 0$.

**2. Symmetric Fuzzy Relation** : A fuzzy relation $R$ on $A$ is said to be *symmetric* iff $\mu_R(x, y) = \mu_R(y, x) \quad \forall x, y \in U$.

**3. Transitive Fuzzy Relation** : A fuzzy relation $R$ on $A$ is said to be *transitive* iff $R \supseteq R \circ R$ where $R \circ R$ is defined by the membership function

$$\mu_{R \circ R}(x, y) = \max_{u \in U} \min \{ \mu_R(x, u), \mu_R(u, y) \}$$

**Definition 2.2** *Let $A = (U, R)$ be an approximation space and $X (\subseteq U)$ be a rough set in $A$. Let $Y$ be a fuzzy set in $U$ with membership grade $\mu_Y : U \to [0, 1]$. Then $Y$ is said to be fuzzy rough set in $A$ if the following*

*conditions holds :*

*(i)* $\mu_Y(x) = 1 \quad \forall\, x \in \underline{X}$

*(ii)* $\mu_Y(x) = 0 \quad \forall\, x \in (U \setminus \overline{X})$

*(iii)* $0 < \mu_Y(x) < 1 \quad \forall\, x \in (\overline{X} \setminus \underline{X})$

## 3 Fuzzy Rough Relation

**Definition 3.1** *Let $A = (U, R)$ be an approximation space and $X\,(\subseteq U)$ be a rough set in A. Let Y be a fuzzy set in A characterized by the membership function $\mu_Y : U \to [0, 1]$. Then the product $Y \times Y$ is defined by the membership function*

$$\mu_{Y \times Y}(x, y) = \min\{\mu_Y(x), \mu_Y(y)\}$$

$$\forall\,(x, y) \in U \times U.$$

A fuzzy rough relation $R$ on $Y$ is a fuzzy rough subset of $Y \times Y$,

i.e. $\mu_R(x, y) \leq \mu_{Y \times Y}(x, y)$

$$\forall\,(x, y) \in U \times U.$$

Satisfying the followings :

(i) $\mu_R(x, y) = 1 \quad \forall\,(x, y) \in \underline{X \times X}$

(ii) $\mu_R(x, y) = 0 \quad \forall\,(x, y) \in [U \times U \setminus \overline{X \times X}]$

(iii) $0 < \mu_R(x, y) < 1$

$$\forall\,(x, y) \in [\overline{X \times X} \setminus \underline{X \times X}]$$

where

$\underline{X \times X} = \{(x, y) \in U \times U : [x, y]_S \subseteq X\}$ and

$\overline{X \times X} = \{(x, y) \in U \times U : [x, y]_S \cap X \neq \varphi\}$.

**Proposition 3.2** *Let $\mu_{R_1}, \mu_{R_2}$ be two FR relations on $X\,(\subseteq U)$ and $\mu$ be pre-defined membership grade on U. Then $\mu_{R_1} \wedge \mu_{R_2}$ is so also.*

**Proof**: Let $\mu' = \mu_{R_1} \wedge \mu_{R_2}$. Then

$$\mu'(x, y) = \left(\mu_{R_1} \wedge \mu_{R_2}\right)(x, y)$$

$$= \min\{\mu_{R_1}(x, y), \mu_{R_2}(x, y)\}.$$

Since $\mu_{R_1}, \mu_{R_2}$ are two FR relations,

$$\mu_{R_1}(x, y) = 1 = \mu_{R_2}(x, y)$$

$$\forall\,(x, y) \in \underline{X \times X}.$$

Therefore

$$\min\{\mu_{R_1}(x, y), \mu_{R_2}(x, y)\} = 1$$

$$\forall\,(x, y) \in \underline{X \times X}$$

$$\Rightarrow \mu'(x, y) = 1 \quad \forall\,(x, y) \in \underline{X \times X}.$$

Also, since $\mu_{R_1}, \mu_{R_2}$ are two FR relations,

$$\mu_{R_1}(x, y) = 0 = \mu_{R_2}(x, y)$$

$$\forall\,(x, y) \in [U \times U \setminus \overline{X \times X}]$$

So,

$$\mu'(x, y) = 0 \quad \forall\,(x, y) \in [U \times U \setminus \overline{X \times X}].$$

Again, since

$$0 < \mu_{R_1}(x, y), \mu_{R_2}(x, y) < 1$$

$$\forall\,(x, y) \in [\overline{X \times X} \setminus \underline{X \times X}]$$

$$\Rightarrow 0 < \min\{\mu_{R_1}(x, y), \mu_{R_2}(x, y)\} < 1$$

$$\forall\,(x, y) \in [\overline{X \times X} \setminus \underline{X \times X}]$$

i.e., $0 < \mu'(x, y) < 1$

$$\forall\,(x, y) \in [\overline{X \times X} \setminus \underline{X \times X}].$$

Hence $\mu'$ is a FR relation on $X$.

**Proposition 3.3** *Let $\mu_{R_1}, \mu_{R_2}$ be two FR relations on $X\,(\subseteq U)$ and $\mu$ is pre-defined membership grade on U. Then $\mu_{R_1} \vee \mu_{R_2}$ is so also.*

**Proof**: Let $\mu' = \mu_{R_1} \vee \mu_{R_2}$.

Then $\mu'(x, y) = (\mu_{R_1} \vee \mu_{R_2})(x, y)$
$= \max\{\mu_{R_1}(x, y), \mu_{R_2}(x, y)\}$.

Since $\mu_{R_1}, \mu_{R_2}$ are two FR relations,

$\mu_{R_1}(x, y) = 1 = \mu_{R_2}(x, y)$
$\forall (x, y) \in \underline{X \times X}$.

Therefore

$\max\{\mu_{R_1}(x, y), \mu_{R_2}(x, y)\} = 1$
$\forall (x, y) \in \underline{X \times X}$

$\Rightarrow \mu'(x, y) = 1 \quad \forall (x, y) \in \underline{X \times X}$.

Also, since $\mu_{R_1}, \mu_{R_2}$ are two FR relations,

$\mu_{R_1}(x, y) = 0 = \mu_{R_2}(x, y)$
$\forall (x, y) \in [U \times U \setminus \overline{X \times X}]$

So,

$\mu'(x, y) = 0 \quad \forall (x, y) \in [U \times U \setminus \overline{X \times X}]$

Again, since
$0 < \mu_{R_1}(x, y), \mu_{R_2}(x, y) < 1$
$\forall (x, y) \in [\overline{X \times X} \setminus \underline{X \times X}]$

$\Rightarrow 0 < \max\{\mu_{R_1}(x, y), \mu_{R_2}(x, y)\} < 1$
$\forall (x, y) \in [\overline{X \times X} \setminus \underline{X \times X}]$

i.e., $0 < \mu'(x, y) < 1$
$\forall (x, y) \in [\overline{X \times X} \setminus \underline{X \times X}]$.

Hence $\mu'$ is a FR relation on $X$.

**Proposition 3.4 ( Algebraic Product )**: Let $\mu_{R_1}, \mu_{R_2}$ be two FR relation on $X (\subseteq U)$ and $\mu$ is pre-defined membership grade on $U$. Then $\mu_{R_1} \cdot \mu_{R_2}$ is FR relation, where $\mu_{R_1} \cdot \mu_{R_2}$ is defined by

$(\mu_{R_1} \cdot \mu_{R_2})(x, y)$
$= \mu_{R_1}(x, y) \mu_{R_2}(x, y)$
$\forall (x, y) \in U \times U$.

**Proof**: Obvious.

**Proposition 3.5 ( Algebraic Sum )**: Let $\mu_{R_1}, \mu_{R_2}$ be two FR relation on $X (\subseteq U)$ and $\mu$ is pre-defined membership grade on $U$. Then $\mu_{R_1} \oplus \mu_{R_2}$ is FR relation, where $\mu_{R_1} \oplus \mu_{R_1}$ is defined by

$(\mu_{R_1} \oplus \mu_{R_1})(x, y) = \mu_{R_1}(x, y) + \mu_{R_2}(x, y)$
$- \mu_{R_1}(x, y) \mu_{R_2}(x, y), \forall (x, y) \in U \times U$

**Proof**: Obvious.

**Definition 3.6** An FR relation $\mu_{R_1}$ on $X (\subseteq U)$ is said to be reflexive FR relation if $\mu_{R_1}(x, y) = 1 \quad \forall x \in U$ such that $\mu(x) > 0$ where $\mu : U \to [0, 1]$ is pre-defined membership function on $U$.

**Definition 3.7** An FR relation $\mu_{R_1}$ on $X (\subseteq U)$ is said to be reflexive FR relation of order $\alpha$ if $\mu_{R_1}(x, y) = \alpha \quad \forall x \in U$ such that $\mu(x) > 0$ where $\mu : U \to [0, 1]$ is pre-defined membership function on $U$.

**Definition 3.8** An FR relation $\mu_{R_1}$ on $X (\subseteq U)$ is said to be weakly reflexive FR relation if $\mu_{R_1}(x, x) \geq \mu_{R_1}(x, y) \forall y \in U$ and $\mu(x) > 0$ where $\mu : U \to [0, 1]$ is pre-defined membership function on $U$.

**Definition 3.9** An FR relation $\mu_{R_1}$ on

$X (\subseteq U)$ is said to be $w$ – reflexive FR relation if $\mu_{R_1}(x, x) \geq \mu(x) \; \forall \, x \in U$ where $\mu : U \to [0, 1]$ is pre-defined membership function on $U$.

**Proposition 3.10** Let $\mu_{R_1}$, $\mu_{R_2}$ be two reflexive FR relation on $X (\subseteq U)$ and $\mu$ is pre-defined membership grade on $U$. Then $\mu_{R_1} \wedge \mu_{R_2}$, $\mu_{R_1} \vee \mu_{R_2}$ are so also.

**Proof**: It is enough to show that $\mu_{R_1} \wedge \mu_{R_2}$, $\mu_{R_1} \vee \mu_{R_2}$ are reflexive.

Let $\mu' = \mu_{R_1} \wedge \mu_{R_2}$. Since $\mu_{R_1}$, $\mu_{R_2}$ are reflexive on $X$,

$$\mu_{R_1}(x, x) = 1 = \mu_{R_2}(x, x) \; \forall \, x$$

with $\mu(x) > 0$. So,

$$\mu'(x, x) = \min\{\mu_{R_1}(x, x), \mu_{R_2}(x, x)\} = 1$$

$\forall \, x$ with $\mu(x) > 0$.

Thus, $\mu'$ is reflexive on $X$.

Similarly, $\left(\mu_{R_1} \vee \mu_{R_2}\right)(x, x)$

$= \max\{\mu_{R_1}(x, x), \mu_{R_2}(x, x)\} = 1$

$\forall x$ with $\mu(x) > 0$.

Hence, $\mu_{R_1} \vee \mu_{R_2}$ and $\mu_{R_1} \wedge \mu_{R_2}$ are reflexive FR relations.

**Proposition 3.11** ( **Algebraic Product** ): Let $\mu_{R_1}, \mu_{R_2}$ be two reflexive FR relation on $X (\subseteq U)$ and $\mu$ is pre-defined membership grade on $U$. Then $\mu_{R_1} \cdot \mu_{R_2}$ is reflexive FR relation, where $\mu_{R_1} \cdot \mu_{R_2}$ is defined by

$$\left(\mu_{R_1} \cdot \mu_{R_2}\right)(x, y) = \mu_{R_1}(x, y) \, \mu_{R_2}(x, y)$$

$\forall (x, y) \in U \times U$.

**Proof**: Obvious.

**Proposition 3.12** ( **Algebraic Sum** ): Let $\mu_{R_1}$, $\mu_{R_2}$ be two reflexive FR relation on $X (\subseteq U)$ and $\mu$ is pre-defined membership grade on $U$. Then $\mu_{R_1} \oplus \mu_{R_2}$ is reflexive FR relation, where $\mu_{R_1} \oplus \mu_{R_2}$ is defined by

$$\left(\mu_{R_1} \oplus \mu_{R_2}\right)(x, y) = \mu_{R_1}(x, y) + \mu_{R_2}(x, y) - \mu_{R_1}(x, y) \, \mu_{R_2}(x, y), \; \forall (x, y) \in U \times U.$$

**Proof**: Obvious.

## 4 Composition of two Fuzzy Rough relations

Let $\mu_{R_1}$, $\mu_{R_2}$ be two FR relations on $X (\subseteq U)$. Then composition of these two FR relations, denoted by $\mu_{R_1} \circ \mu_{R_2}$, is defined by

$$\left(\mu_{R_1} \circ \mu_{R_2}\right)(x, y)$$

$$= \max_{u \in U} \min\{\mu_{R_1}(x, u), \mu_{R_2}(u, y)\}$$

$\forall \, x, y \in U$.

**Proposition 4.1** Composition of two FR relations on $X (\subseteq U)$ is associative.

**Proof**: Obvious.

**Proposition 4.2** Let $\mu_{R_1}$, $\mu_{R_2}$ be two FR relations on $X (\subseteq U)$ and $\mu$ be pre-defined membership grade on $U$. Then $\mu_{R_1} \circ \mu_{R_2}$ is FR relation.

**Proof**: Since $\mu_{R_1}$, $\mu_{R_2}$ are two FR relations

$$\mu_{R_1}(x, y) = 1 = \mu_{R_2}(x, y)$$

$\forall \, (x, y) \in \underline{X \times X}.$

Let $(x, y) \in \underline{X \times X}$. Now,

$\left( \mu_{R_1} \circ \mu_{R_2} \right)(x, y)$

$= \max_{u \in U} \min \left\{ \mu_{R_1}(x, u), \mu_{R_2}(u, y) \right\}$

$= 1$. This holds for all $(x, y) \in \underline{X \times X}$.

Let $(x, y) \in \left[ U \times U \setminus \overline{X \times X} \right]$. So,

$\mu_{R_1}(x, y) = 0 = \mu_{R_2}(x, y)$ Then

$\left( \mu_{R_1} \circ \mu_{R_2} \right)(x, y)$

$= \max_{u \in U} \min \left\{ \mu_{R_1}(x, u), \mu_{R_2}(u, y) \right\} = 0$

Again, since $0 < \mu_{R_1}(x, y), \mu_{R_2}(x, y) < 1$

$\forall (x, y) \in [ \overline{X \times X} \setminus \underline{X \times X} ]$,

$\Rightarrow 0 < \max_{u \in U} \min \left\{ \mu_{R_1}(x, u), \mu_{R_2}(u, y) \right\} < 1$

i.e. $0 < \left( \mu_{R_1} \circ \mu_{R_2} \right)(x, y) < 1$

$\forall (x, y) \in [ \overline{X \times X} \setminus \underline{X \times X} ]$.

Hence, $\mu_{R_1} \circ \mu_{R_2}$ is also FR relation on $X$.

**Proposition 4.3** Let $\mu_{R_1}, \mu_{R_2}$ be two reflexive FR relations on $X (\subseteq U)$ and $\mu$ be pre-defined membership grade on $U$. Then $\mu_{R_1} \circ \mu_{R_2}$ is also reflexive FR relation.

**Proof**: By proposition 4.1, $\mu_{R_1} \circ \mu_{R_2}$ is a FR relation on $X$. So, it suffices to show that $\mu_{R_1} \circ \mu_{R_2}$ is a reflexive FR relation on $X$.

We see that,

$\left( \mu_{R_1} \circ \mu_{R_2} \right)(x, x)$

$= \max_{u \in U} \min \left\{ \mu_{R_1}(x, u), \mu_{R_2}(u, x) \right\}$

$= \max_u \left\{ \max_{u \neq x} \min \left\{ \mu_{R_1}(x, u), \mu_{R_2}(u, x) \right\}, \right.$

$\left. \max_{u = x} \min \left\{ \mu_{R_1}(x, u), \mu_{R_2}(u, x) \right\} \right\}$

$= 1$ [ since $\mu_{R_1}, \mu_{R_2}$ are reflexive FR relations ]

Hence $\mu_{R_1} \circ \mu_{R_2}$ is reflexive FR relation.

**Definition 4.4** An FR relation $\mu_{R_1}$ on $X (\subseteq U)$ is said to be symmetric FR relation if $\mu_{R_1}(x, y) = \mu_1(y, x) \ \forall \ x, y \in U$.

**Definition 4.5** An FR relation $\mu_{R_1}$ on $X (\subseteq U)$ is said to be transitive FR relation if $\mu_{R_1} \supseteq \mu_{R_1} \circ \mu_{R_1}$, where '$\circ$' is defined by

$\left( \mu_{R_1} \circ \mu_{R_1} \right)(x, y)$

$= \max_{u \in U} \min \left\{ \mu_{R_1}(x, u), \mu_{R_1}(u, y) \right\}$

$\forall \ x, y \in U$.

**Proposition 4.6** Let $\mu_{R_1}$ be symmetric and transitive FR relation then $\mu_{R_1}$ is weakly reflexive FR relation.

**Proof**: Since $\mu_{R_1}$ is symmetric,

$\mu_{R_1}(x, y) = \mu_{R_1}(y, x) \ \forall \ x, y \in U$.

Again, since $\mu_{R_1}$ is transitive,

$\mu_{R_1}(x, y) \geq \max_{u \in U} \min \left\{ \mu_{R_1}(x, u), \mu_{R_1}(u, y) \right\}$

Taking $y = x$, we get

$\mu_{R_1}(x, x) \geq \max_{u \in U} \min \left\{ \mu_{R_1}(x, u), \mu_{R_1}(u, x) \right\}$

$= \max_{u \in U} \min \left\{ \mu_{R_1}(x, u), \mu_{R_1}(x, u) \right\}$

$= \max_{u \in U} \mu_{R_1}(x, u)$

Therefore, $\mu_{R_1}(x,x) \geq \mu_{R_1}(x,y)$ for all $y \in U$ where $\mu(x) > 0$.

This shows that $\mu_{R_1}$ is weakly reflexive FR relation on $X$.

**Proposition 4.7** Let $\mu_{R_1}$, $\mu_{R_2}$ be w-reflexive FR relations then $\mu_{R_1} \vee \mu_{R_2} \subseteq \mu_{R_1} \circ \mu_{R_2}$.

**Proof**: We see that,

$\left(\mu_{R_1} \circ \mu_{R_1}\right)(x,y)$

$= \max_{u \in U} \min\left\{\mu_{R_1}(x,u), \mu_{R_1}(u,y)\right\}$

$\geq \min\left\{\mu_{R_1}(x,x), \mu_{R_2}(x,y)\right\}$

$\geq \min\left\{\mu(x), \mu_{R_2}(x,y)\right\}$

[ as $\mu_{R_1}$ is w – reflexive ]

Again,

$\mu_{R_2}(x,y) \leq \min\{\mu(x), \mu(y)\} \leq \mu(x)$

So, $\left(\mu_{R_1} \circ \mu_{R_2}\right)(x,y) \geq \mu_{R_2}(x,y)$.

Therefore, $\mu_{R_2} \subseteq \mu_{R_1} \circ \mu_{R_2}$ Similarly, we can show that $\mu_{R_1} \subseteq \mu_{R_1} \circ \mu_{R_2}$.

Hence, $\mu_{R_1} \vee \mu_{R_2} \subseteq \mu_{R_1} \circ \mu_{R_2}$.

**Definition 4.8** *An FR relation is on $X (\subseteq U)$ is said to be a **similitude** FR relation if it is reflexiv, symmetric and transitive.*

**Theorem 4.9** Let $X$ be a fuzzy rough subset of universe $U$ where $\mu$ be a membership function on $U$. Also let $\mu_{R_1}$, a similitude FR relation of order $\alpha$. Then for each $x \in U$, with $\mu(x) > 0$, $\exists$ a fuzzy rough subset of $X$ determined by membership function $\mu_{R_x}$ satisfying the followings :

(i) $\mu_{R_x}(x) = \alpha$

(ii) $\mu_{R_x}(y) = \mu_{R_y}(x)$

(iii) $\mu_{R_x}(y) > 0$, $\mu_{R_y}(z) > 0$

$\Rightarrow \mu_{R_x}(z) > 0$

(iv) $\mu_{R_x}(u) = 0 \Rightarrow \mu_{R_x} \wedge \mu_{R_y} = \varphi$.

**Proof**: We see that $\mu_{R_1}(x,y) \leq \min\{\mu(x), \mu(y)\} \quad \forall x, y \in U$

For each $x \in U$, with $\mu(x) > 0$, we defined

$\mu_{R_x}(y) = \mu_{R_1}(x,y) \quad \forall y \in U$. We note that, $\mu_{R_x}(y) = \mu_{R_1}(x,y)$

$\leq \min\{\mu(x), \mu(y)\} \leq \mu(y) \quad \forall y \in U$.

i.e. $\mu_{R_x}(y) \leq \mu(y) \quad \forall y \in U$.

Thus the fuzzy rough set determined by $\mu_{R_x}$ is a fuzzy rough subset of $X$.

(i) $\mu_{R_x}(x) = \mu_{R_1}(x,x) = \alpha$ [ since $\mu_{R_1}$ is similitude of order $\alpha$ ]

(ii) $\mu_{R_x}(y) = \mu_{R_x}(x,y)$

$= \mu_{R_x}(y,x) = \mu_{R_y}(x)$

[ as $\mu_{R_1}$ is symmetric ]

(iii) We suppose that $\mu_{R_x}(y) > 0$ and $\mu_{R_y}(z) > 0$

i.e. $\mu_{R_1}(x,y)$, $\mu_{R_1}(y,z) > 0$ \qquad (1)

By transitivity of $\mu_{R_1}$,

$\mu_{R_1}(x,z) \geq \max_{u \in U} \min\left\{\mu_{R_1}(x,u), \mu_{R_1}(u,z)\right\}$

Claim:

$\max_{u \in U} \min\left\{\mu_{R_1}(x,u), \mu_{R_1}(u,z)\right\} > 0$.

For $u = y \, (\in U)$, it becomes

$$\max_{u} \left\{ \max_{u \neq y} \min \left\{ \mu_{R_1}(x, u), \mu_{R_1}(u, x) \right\}, \right.$$

$$\left. \max_{u = y} \min \left\{ \mu_{R_1}(x, u), \mu_{R_1}(u, x) \right\} \right\} > 0$$

[ by (1) ]

So, our claim is justified.

Thus, $\mu_{R_x}(z) = \mu_{R_1}(x, z) > 0$.

(iv) Let $\mu_{R_x}(u) = 0$. We need to show that $\mu_{R_x} \wedge \mu_{R_y} = \varphi$. If possible, let $\left( \mu_{R_x} \wedge \mu_{R_y} \right)(z) > 0$

$\Rightarrow \min \left\{ \mu_{R_x}(z), \mu_{R_y}(z) \right\} > 0$

So, $\mu_{R_x}(z), \mu_{R_y}(z) > 0$.

i.e. $\mu_{R_1}(x, z), \mu_{R_1}(y, z) > 0$. Now,

$\mu_{R_1}(x, u) \geq \max_{t \in U} \min \left\{ \mu_{R_1}(x, t), \mu_{R_1}(t, u) \right\}$.

Taking $t = z$ we have,

$\mu_{R_1}(x, u) = \mu_{R_x}(u) > 0$,

This contradicts $\mu_{R_x}(u) = 0$.

Thus, $\mu_{R_x} \wedge \mu_{R_y} = \varphi$.

This completes the proof.